\documentclass[12pt]{amsart}
\usepackage{amsfonts}
\usepackage{amsfonts, mathdots}
\usepackage{amssymb,latexsym,color}
\usepackage{enumerate}\usepackage{arydshln}
\usepackage{hyperref}% webpage
\makeatletter
\@namedef{subjclassname@2010}{%
  \textup{2010} Mathematics Subject Classification}
\makeatother

\newtheorem{thm}{Theorem}[section]

\theoremstyle{definition}

\numberwithin{equation}{section}

\begin{document}

\baselineskip=17pt

\title[power majorization]{Power majorization between the roots of two polynomials}

\author[M. Lin]{Minghua Lin}
\address{Department of Mathematics \\
  Shanghai University\\ Baoshan district, Shanghai, 200444, China}
\email{m\_lin@i.shu.edu.cn \\ URL: http://my.shu.edu.cn/en/mlin}
\author[G. Sinnamon]{Gord Sinnamon}
\address{Department of Mathematics \\
  University of Western Ontario\\ London, ON, Canada}
\email{sinnamon@uwo.ca\\ URL: http://www-home.math.uwo.ca/$\sim$sinnamon/}

\begin{abstract} It is shown that if two hyperbolic polynomials have a particular factorization into quadratics, then their roots satisfy a power majorization relation whenever key coefficients in their factorizations satisfy a corresponding majorization relation. In particular, a numerical observation by Kleme\v{s}  \cite{Kle11} is confirmed.  
   \end{abstract}



\maketitle

\section{Introduction}
 \label{s1}

Let $x=(x_1, \ldots, x_n)$ and $y=(y_1, \ldots, y_n)$ be two  $n$-tuples of real numbers.  We recall (\cite[Chapter 1]{MOA11}) that $x$ is said to be majorized by $y$ if the sum of $k$   largest entries of $x$ is less or equal to   the sum of $k$  largest entries of $y$, where $k$ ranges from $1$ to $n-1$, and equality holds at $k=n$. 
 There are several useful characterizations of majorization, perhaps the most famous one is due to Hardy, Littlewood and P\'{o}lya \cite[p. 156]{MOA11} which says that  $x$ is   majorized by $y$ if and only if
\begin{eqnarray*}
 \sum_{i=1}^n\varphi(x_i)\le \sum_{i=1}^n\varphi(y_i)
\end{eqnarray*}
for every convex function $\varphi$.

Mainly motivated by the study of $l^p$-means, for two $n$-tuples of positive real numbers, we say that $x$ is power majorized by $y$
 provided that
  \begin{eqnarray*}
   \sum_{i=1}^n x_i^p\le \sum_{i=1}^n y_i^p
  \end{eqnarray*}
 whenever $p\ge 1$, with reversal of the inequality sign when $0<p<1$. 
 
 Clearly, majorization implies power majorization. But the converse is not true. The first example illustrating the difference between   majorization and power majorization was given in \cite{Ben86}. Unlike the rich theory on the majorization relation, little is known about  power majorization.  Indeed, it is in general difficult to determine whether one vector is power majorized by another.  Our investigation in this paper stems from  the following example. 
 
 Suppose that one is interested in
 comparing the $l^p$ norms of the eigenvalues $x =(x_1, \ldots, x_4)$ and $y=(y_1, \ldots,  y_4)$
 respectively of the $4\times 4$ matrices $X$ and $Y$ defined by $X = AA^T$, $Y = BB^T$,
 where \begin{eqnarray*} A=\begin{bmatrix}   1 & 1& 1& 1\\
  0& 1& 1& 0\\
  0& 0& 1 & 0\\
  0& 0&   1& 1
  \end{bmatrix}, \quad B=\begin{bmatrix}  
  1& 1& 1& 1\\
  0& 1& 1& 0\\
  0& 0& 1& 0\\
  0& 0& 1& -1 \end{bmatrix}.
  \end{eqnarray*} 
 On numerical evidence, it was suggested in  \cite{Kle11}  that   $y$ is power majorized by  $x$.  
 The author of \cite{Kle11}  asked for an 
  “enlightening” proof or disproof.  We   confirm this numerical guess.

 A simple calculation gives 
  \begin{eqnarray*} X=\begin{bmatrix} 4  &   2  &   1 &    2\\
         2 &    2   &  1 &    1\\
         1 &    1    & 1   &  1\\
         2 &    1   &  1   &  2
    \end{bmatrix}, \quad Y=\begin{bmatrix} 4  &   2  &   1 &    0\\
         2 &    2   &  1 &    1\\
         1 &    1    & 1   &  1\\
        0 &    1   &  1   &  2
    \end{bmatrix}
    \end{eqnarray*}
    and their characteristic polyonmials
    \begin{eqnarray*} P(t)=(t^2-7t+1)(t^2-2t+1), \\ Q(t)=(t^2-6t+1)(t^2-3t+1). \end{eqnarray*}
  
  It can be shown that in this example $y$ is not majorized by $x$ (as the sum of the two largest   entries of $y$ is larger than the sum of the two largest entries of $x$). We observe that in the above two factors of  $P(t)$ and $Q(t)$, the coefficients of  $t$ satisfy a majorization relation. More precisely, $(6, 3)$ is majorized by $(7, 2)$.   This simple observation hints at the statement of our main result.

   \section{Main Result}
   
 A hyperbolic polynomial is a polynomial whose roots are all real. We refer the interested reader to the classical text  \cite{Obr63}  on the roots of hyperbolic polynomials.  Our main result is the following theorem, which gives a condition for  power majorization between the roots of two hyperbolic polynomials  having a certain factorization.   We remark that relevant studies on the majorization relation between the roots of hyperbolic polynomials  (the so called spectral order) are given in 
 \cite{BS03, Bor06, CPK11, Pei05}.
 
   \begin{thm}\label{main} Let $u=(u_1, \ldots, u_n)$,  $v=(v_1, \ldots, v_n)$, with $u_i, v_i\ge 1$, $i=1, \ldots, n$. Consider two polynomials 
     \begin{eqnarray*} P(t)=\prod_{i=1}^n(t^2-2u_it+1), \\ Q(t)=\prod_{i=1}^n(t^2-2v_it+1).  \end{eqnarray*}
    Let $x =(x_1, \ldots, x_{2n})$ and $y=(y_1, \ldots,  y_{2n})$ be the vectors of the roots of $P(t)$ and $Q(t)$, respectively. If $v$ is majorized by $u$, then $y$   is power majorized by $x$.  
   \end{thm}
   
 We need some basic facts about   Schur-convex functions (\cite[p. 80]{MOA11}). 
    A real-valued function $\varphi$ defined on a set $\mathcal{A}\subset \mathbb{R}^n$
   is said to be Schur-convex on $\mathcal{A}$ if for every $x, y\in \mathcal{A}$, $y$ being majorized by $x$  implies  $\varphi(x) \ge  \varphi(y)$. We say $\varphi$ is Schur-concave if $-\varphi$ is Schur-convex.
   
   Let $\mathcal{I}\subset \mathbb{R}$ be an
   open interval and let $\varphi: \mathcal{I}^n \to  \mathbb{R}$ be continuously differentiable.
   The well-known   Schur  condition (\cite[p. 84]{MOA11}) says that $\varphi$  is Schur-convex on  $\mathcal{I}^n$  if and only if  $\varphi$ is symmetric on $\mathcal{I}^n$ and  $$(x_1-x_2)\left(\frac{\partial \varphi}{\partial x_1}-\frac{\partial \varphi}{\partial x_2}\right)\ge 0. $$
 
 \noindent{\it Proof of Theorem \ref{main}.}  The roots of $P(t)$ are $u_i\pm \sqrt{u_i^2-1}$, $i=1, \ldots, n$ and the roots of $Q(t)$ are $v_i\pm \sqrt{v_i^2-1}$, $i=1, \ldots, n$. 
  We need to show that for any fixed $p\ge 1$, 
  \begin{eqnarray*}
  	&&  \sum_{i=1}^n \left(\left(u_i+\sqrt{u_i^2-1}\right)^p+\left(u_i-\sqrt{u_i^2-1}\right)^p\right)\\ && \qquad \ge  \sum_{i=1}^n \left(\left(v_i+\sqrt{v_i^2-1}\right)^p+\left(v_i-\sqrt{v_i^2-1}\right)^p\right) 
  \end{eqnarray*} and that the inequality reverses for $0<p<1$. 
  
  This would follow if one could show
  \begin{eqnarray*}  \varphi(u):= \sum_{i=1}^n \left(\left(u_i+\sqrt{u_i^2-1}\right)^p+\left(u_i-\sqrt{u_i^2-1}\right)^p\right),  \end{eqnarray*}  where $u\in [1, \infty)^n$, is Schur convex for $p\ge 1$ and is Schur concave for $0<p<1$.   
  
Assume first that $p\ge 1$.  Clearly,   $\varphi(u)$ is symmetric; by   Schur's  condition, it remains to show that 
  $$(u_1-u_2)\left(\frac{\partial \varphi}{\partial u_1}-\frac{\partial \varphi}{\partial u_2}\right)\ge 0. $$
  
  Without loss of generality, we assume $u_1\ge u_2$. Then it suffices to show  \begin{eqnarray}\label{e21}\frac{\partial \varphi}{\partial u_1}-\frac{\partial \varphi}{\partial u_2}\ge 0. 
  \end{eqnarray}
  Compute 
  \begin{eqnarray*}  
  	\frac{\partial \varphi}{\partial u_1}-\frac{\partial \varphi}{\partial u_2}
  	&=& p\frac{(u_1+\sqrt{u_1^2-1})^p-(u_1-\sqrt{u_1^2-1})^p}{\sqrt{u_1^2-1}} \\&& \qquad  -p\frac{(u_2+\sqrt{u_2^2-1})^p-(u_2-\sqrt{u_2^2-1})^p}{\sqrt{u_2^2-1}}.\end{eqnarray*}
  Thus, (\ref{e21}) would follow if we can show 
  \begin{eqnarray*}  
  	g(t):=\frac{(t+\sqrt{t^2-1})^p-(t-\sqrt{t^2-1})^p}{\sqrt{t^2-1}}  \end{eqnarray*} is an increasing function for $t>1$. 
  
  A simple calculation gives  \begin{eqnarray*}  
  	g'(t)=\frac{p(t+\sqrt{t^2-1})^p+p(t-\sqrt{t^2-1})^p-\Big((t+\sqrt{t^2-1})^p-(t-\sqrt{t^2-1})^p\Big)\frac{t}{\sqrt{t^2-1}}}{t^2-1}.\end{eqnarray*} 
  With  $\theta=\sqrt{t^2-1}/t$, this becomes
     \begin{eqnarray*}  \frac{t^p}{t^2-1}\Big(p(1+\theta)^p+p(1-\theta)^p-((1+\theta)^p-(1-\theta)^p)/\theta\Big).
     	\end{eqnarray*} 
 
 To see that $g'(t)\ge 0$, it suffices to show that for $0<\theta<1$,  \begin{eqnarray*}  (p\theta-1)(1+\theta)^p+(p\theta+1)(1-\theta)^p\ge 0.
 \end{eqnarray*} 
When $p\theta\ge 1$ this is clear. When $p\theta<1$ it is equivalent to showing 
 \begin{eqnarray*}  
 	h(\theta):=\ln((p\theta+1)(1-\theta)^p)-\ln((1-p\theta)(1+\theta)^p)\ge 0.   \end{eqnarray*} But 
 $h(0)=0$ and 
 \begin{eqnarray*}  
 	h'(\theta)=\frac{2p}{1-p^2\theta^2}-\frac{2p}{1-\theta^2}\ge 0   \end{eqnarray*}
so $h(\theta)\ge 0$ for $0<\theta<1/p$. This completes the proof of the $p\ge 1$ case. 
 
 If $0<p<1$, the argument is similar to the proof of the $p\ge 1$ case. It suffices to show  $g(t)$ defined above is decreasing for $t>1$. To see that $g'(t)\le 0$, we need to verify that for $0<\theta<1$,  \begin{eqnarray*}  (p\theta-1)(1+\theta)^p+(p\theta+1)(1-\theta)^p\le 0.
 \end{eqnarray*} Equivalently, $h(\theta)$ defined above should be nonpositive. But in this case, $h(0)=0$ and   \begin{eqnarray*}  
 h'(\theta)=\frac{2p}{1-p^2\theta^2}-\frac{2p}{1-\theta^2}\le 0.   \end{eqnarray*}
 So the proof of Theorem \ref{main} is complete. \qed

\subsection*{Acknowledgments} {\small The work of M. Lin is supported by the National Natural Science Foundation of China. The work of G. Sinnamon is supported  by the Natural Sciences and Engineering Research Council of Canada. }

\end{document}